\documentclass[conference,letterpaper]{IEEEtran}
\IEEEoverridecommandlockouts
\usepackage{cite}
\usepackage{amsmath,amssymb,amsfonts}
\usepackage{algorithmic}
\usepackage{graphicx}
\usepackage{textcomp}
\usepackage{xcolor}
\usepackage{array}
\usepackage{hyperref}
\usepackage{parskip}

\hypersetup{
    colorlinks=true,
    linkcolor=black,
    citecolor=black,
    urlcolor=blue
}

\newcommand{\numPipelines}{P}              
\newcommand{\numCandidateStations}{Q}      
\newcommand{\numChosenStations}{S}         
\newcommand{\numRobots}{R}                 

\newcommand{\transitSpeed}{\nu_\textrm{tra}}        
\newcommand{\inspectSpeed}{\nu_\textrm{ins}}        

\newcommand{\numSamples}[1]{N_{#1}}        

\newcommand{\entryCoord}[1]{\mathbf{b}_{#1}}      

\newcommand{\anomalyCoord}{\mathbf{a}_{p, i}}

\newcommand{\stationCoord}[1]{\mathbf{q}_{#1}}    

\def\tMax{T_{\max}}             
\def\tAvg{T_{\mathrm{avg}}}     

\def\epsParam{\varepsilon}         
\begin{document}

\begin{titlepage}
\thispagestyle{empty}

\centering

\vspace*{12mm}

{\Large\bfseries Author Accepted Manuscript\par}

\vspace{18mm}

{\LARGE\bfseries
\begin{minipage}{0.96\textwidth}
\centering
Optimal Placement of Docking Stations and Resident AUVs for\\
Subsea Pipeline Inspection
\end{minipage}
\par}

\vspace{16mm}

{\large
Gabrielė Kasparavičiūtė\\[2mm]
Pasquale Grippa\\[2mm]
Kjetil Skaugset\\[2mm]
Asgeir J. Sørensen\\[2mm]
Martin Ludvigsen
\par}

\vfill

\rule{0.88\textwidth}{0.5pt}

\vspace{7mm}

\begin{minipage}{0.88\textwidth}
\centering
\small

This document is the peer-reviewed Author Accepted Manuscript of an
article accepted for publication in the
\textit{IEEE Journal of Oceanic Engineering}.

\vspace{5mm}

\textbf{Published article}

\vspace{2mm}

G.~Kasparavičiūtė, P.~Grippa, K.~Skaugset,
A.~J.~Sørensen, and M.~Ludvigsen,
``Optimal Placement of Docking Stations and Resident AUVs for
Subsea Pipeline Inspection,''
\textit{IEEE Journal of Oceanic Engineering}, 2026.

\vspace{3mm}

DOI:
\href{https://doi.org/10.1109/JOE.2026.3694005}
{\texttt{10.1109/JOE.2026.3694005}}

\vspace{3mm}

Official version:
\href{https://ieeexplore.ieee.org/document/11570255}
{\texttt{ieeexplore.ieee.org/document/11570255}}

\end{minipage}

\vspace{7mm}

\rule{0.88\textwidth}{0.5pt}

\vspace{7mm}

\begin{minipage}{0.88\textwidth}
\footnotesize

\textbf{\copyright{} 2026 IEEE.}
Personal use of this material is permitted. Permission from IEEE must
be obtained for all other uses, in any current or future media,
including reprinting or republishing this material for advertising or
promotional purposes, creating new collective works, for resale or
redistribution to servers or lists, or reuse of any copyrighted
component of this work in other works.

\end{minipage}

\vspace*{8mm}

\end{titlepage}

\title{Optimal Placement of Docking Stations and Resident AUVs for Subsea Pipeline Inspection
}

\author{
\IEEEauthorblockN{
Gabrielė Kasparavičiūtė,
Pasquale Grippa,
Kjetil Skaugset,
Asgeir J. Sørensen,
Martin Ludvigsen
}

\IEEEauthorblockA{
Department of Marine Technology, Norwegian University of Science and Technology (NTNU), Trondheim, Norway
}

\IEEEauthorblockA{
Pasquale Grippa is with Lakeside Labs GmbH, Klagenfurt, Austria.
}

\IEEEauthorblockA{
Kjetil Skaugset is also with Equinor ASA, Norway.
}

\IEEEauthorblockA{
Corresponding author: gabriele.kasparaviciute@ntnu.no
}
}

\maketitle

\begin{abstract}
A two‐stage mixed-integer linear programming framework is introduced for subsea pipeline incident response planning, jointly optimizing Subsea Docking Plate (SDP) placement and resident autonomous underwater vehicle allocation to minimize both maximum and average response times under spatial uncertainty. Phase 1 minimizes the maximum response time across all potential leak locations. Phase 2 reduces the average response time subject to the maximum bound. A case study based on the Johan Sverdrup oil and gas field in Norway, with 9 SDP options and 33 pipelines, shows that just a few well-placed SDPs are enough to achieve top performance. A cost versus time Pareto analysis reveals the trade-off between deployment expenditure and response efficiency, providing actionable guidance for designing resilient and cost-effective subsea inspection networks. 

\end{abstract}


\section{INTRODUCTION}

Subsea conduit systems, generally known as pipelines or flowlines, are vital for offshore oil and gas production, frequently extending tens to hundreds of kilometers in deep water. This study employs the term ``pipeline'' as a general designation for both types, despite modest variations in industry usage. These conduits transport manufacturing fluids and are essential assets necessitating meticulous monitoring and maintenance.

Maintaining pipeline integrity involves monitoring parameters such as pressure to detect anomalies. Upon detecting a potential anomaly, operators typically shut down the affected pipeline as a precautionary measure, dispatching Remotely Operated Vehicles (ROVs) to inspect and resolve the issue. Such shutdowns incur significant financial losses, emphasizing the importance of minimizing the response time which is the duration from anomaly detection to acquiring conclusive inspection data. 

While ROV operations are effective, they are limited by their dependence on favorable weather conditions, as well as the high costs and greenhouse gas emissions associated with support ships \cite{norskpetro}. Furthermore, their response time is long due to the time required to mobilize the ship and deploy the ROV. Autonomous Underwater Vehicles (AUVs) were introduced in recent years, expanding inspection capability with untethered operation and increased spatial coverage. However, AUVs still require periodic surfacing to recharge batteries and transmit data which introduces delays and operational inefficiencies. Recent research suggests that future maritime operations will involve coordinated robotic organizations made up of heterogeneous autonomous platforms that operate consistently across domains. In this context, subsea docking infrastructure is a crucial component that facilitates constant underwater presence and quick response to infrastructure anomalies \cite{skaugset2025autonomous}.

Resident Autonomous Underwater Vehicles (RAUVs), also referred to as Underwater Intervention Drones (UIDs) have emerged as a promising advancement to address these limitations. RAUVs operate continuously underwater, supported by Subsea Docking Plates (SDPs) that provide infrastructure for recharging and data transmission entirely beneath the sea surface, which eliminates resurfacing \cite{vasilijevic2023infrastructure}. The resident approach significantly improves response times and operational reliability. Optimally placing SDPs and assigning RAUVs to these docking stations is crucial for achieving rapid response times and comprehensive coverage of the pipeline network. 

This paper introduces a unified two-phase Mixed-Integer Linear Program (MILP) tailored to the joint problem of SDP placement and RAUV allocation. In Phase 1, the MILP minimizes the maximum round-trip mission time (approach, inspection, and return) over all discretized points along the pipeline network. Phase 2 then reduces the average mission time without violating the Phase 1 maximum. To capture spatial uncertainty in leak locations, each pipeline is discretized into a set of inspection points, each treated as a potential anomaly site. We validate the approach on a realistic case study of 33 pipelines and 9 candidate docking stations, including cost-effectiveness comparisons across a range of RAUV and SDP deployment choices.  

We focus on a MILP formulation because it guarantees the identification of a global optimum. Achieving the global optimum is  particularly important in our use case, where the deployment costs are high and the response time (performance metric) is critical. Furthermore, the MILP model introduced in this paper establishes an optimal benchmark that can be used to evaluate the performance of  heuristic and metaheuristic approaches, such as Ant Colony Optimization. While these methods generally offer computational advantages for large problem instances, they do not provide guarantees of optimality.

The remainder of the paper is organized as follows. Section~\ref{related_work} surveys related placement and incident response models. Sections ~\ref{problem_statement} and ~\ref{solution_method} formalize the RAUV-SDP assignment problem and present the MILP method. Section~\ref{results} reports on solve times and cost trade-offs. Section~\ref{discussion} discusses the results. The conclusion in Section~\ref{conclusion} provides directions for scaling to even larger subsea infrastructures and future work.


\section{RELATED WORK}
\label{related_work}

\begin{table*}[ht]
\begin{center}
\caption{Concise comparison of related optimisation works and the gap addressed in this paper.}
\label{tab:related_work}
\begin{tabular}{|c|c|>{\centering\arraybackslash}p{3.1cm}|%
                >{\centering\arraybackslash}p{3.1cm}|%
                >{\centering\arraybackslash}p{3.1cm}|c|}
\hline
\textbf{Paper} & \textbf{Ref.} &
\shortstack{\textbf{Primary}\\\textbf{Objective}} &
\shortstack{\textbf{Method}\\(\textit{model})} &
\shortstack{\textbf{Uncertainty}\\\textbf{Model}} &
\shortstack{\textbf{Sub\-sea}\\\textbf{Context}} \\
\hline
Hong et al. (2023) & \cite{hong2023minlp} &
Minimize layout cost of subsea network &
MINLP with geo-heuristic &
Deterministic: fixed wells, obstacles &
\checkmark \\
\hline
Liu et al. (2022) & \cite{liu2022subsea} &
Reduce flowline cost &
Binary LP after clustering &
Deterministic: fixed well locations &
\checkmark \\
\hline
Wang et al. (2022) & \cite{wang2022ems} &
Expected EMS response cost &
2-stage stochastic program &
Scenario-based demand/travel time &
- \\
\hline
Xue et al. (2023) & \cite{xue2023vessel} &
Minimize vessel travel cost with UAV delays &
Fuzzy MILP + genetic–tabu heuristic &
Fuzzy UAV travel times (triangular) &
\checkmark \\
\hline
Li et al. (2025) & \cite{li2025two} &
Facility–cost–impact trade-off &
Robust opt. with C\&CG &
Budgeted uncertainty in demand &
- \\
\hline
\textbf{This paper} &  &
Min/max RAUV mission time &
\textbf{2-phase MILP with anomaly uncertainty} &
Anomaly location unknown &
\checkmark \\
\hline
\end{tabular}
\end{center}
\end{table*}

The problem of selecting a subset of charging stations and assigning resident RAUVs to minimize maximum and average inspection response times lies at the intersection of facility location, vehicle allocation, and incident response optimization. Facility location techniques were used in both terrestrial \cite{wang2022ems, li2025two}, and subsea applications \cite{hong2023minlp}, taking into account context specific features such as uncertainty in anomaly location and the need for both maximum and average response time performance during underwater pipeline inspection (see Table \ref{tab:related_work}). 

Facility placement and vehicle allocation, which are critical components of subsea incident response planning, originate from larger operational planning concepts used in logistics and emergency services. Mixed-integer formulations for joint facility location and vehicle allocation have been used in logistics and emergency services \cite{farahani2010multiple, dukkanci2024facility}. Early covering location and $p$-median models \cite{boujelben2016milp, arabahmadi2023facility} balance the number of open facilities with service level metrics. Covering location (or $p$-center) formulations minimize the maximum distance from a demand point to its nearest facility, while $p$-median formulations minimize the total (or average) distance within a certain facility budget.

Researchers in \cite{onut2008two} propose a two-phase possibilistic programming approach for supplier assessment with fuzzy cost and quality data. In Phase 1, a fuzzy membership operator is used to maximize the minimum satisfaction level across all objectives. In Phase 2, the guaranteed membership level is enforced as a lower bound constraint, and the average membership is optimized using a compensatory operator to improve overall performance. The proposed MILP model in this paper follows a similar two-phase philosophy but uses deterministic travel times. Instead of relying on fuzzy or scenario based proxies, this paper computes the expected inspection and return times for each segment to ensure full realism.

Relatively few studies target underwater incident response. Authors in \cite{davis2015subsea} model acoustic communication constraints in subsea networks, but without joint station placement or multi‐phase maximum/average objectives. In earlier work \cite{kasparavivciute2025energy}, route planning was considered for a single RAUV servicing six pipelines from a fixed docking station. That model lacked SDP placement and vehicle assignment decisions capabilities which are introduced here.

The papers listed in Table \ref{tab:related_work} illustrate various approaches to facility location and vehicle allocation under different contexts and uncertainty models. Certain studies presume deterministic input data and disregard uncertainty, whereas others integrate diverse sources of uncertainty (e.g., such as variable demand, trip time, or inspection duration) into the decision-making framework. Authors of \cite{hong2023minlp} and \cite{liu2022subsea} both address subsea infrastructure (e.g., pipeline network topology, subsea wells, manifolds, processing terminals) optimization and focus primarily on cost rather than response times and do not explicitly consider anomaly locations. Researchers in \cite{xue2023vessel} and \cite{li2025two} incorporate sophisticated uncertainty models (stochastic and robust, respectively) in terrestrial incident contexts without specific subsea considerations. Scientists in \cite{xue2023vessel} address inspection timing within subsea scenarios and employ fuzzy uncertainty models. In contrast, this paper assumes that anomalies are equally likely to occur in any pipeline section. It considers both maximum and average response time optimization within a unified MILP framework that is explicitly tailored to subsea RAUV inspection tasks. 


Furthermore, applying the MILP model to different numbers of SDPs and RAUVs and associating these numbers with acquisition and maintenance costs allow us to link the cost of a system to its response time. Previous study \cite{grippa2018} demonstrates a similar trade-off in drone delivery systems, specifically the balance between infrastructure investment (number of drones and depots) and operational efficiency (average delivery time).


Within these contexts, the proposed two-phase mixed-integer linear programming model constitutes as a novel contribution with the capacity to optimize the placement of SDPs and the allocation of RAUVs with consideration of the uncertainty of the location of anomalies under a subsea inspection scenario.

\section{PROBLEM DESCRIPTION}
\label{problem_statement}
\begin{figure}[ht]
 \centering
 \includegraphics[width=0.75\columnwidth]{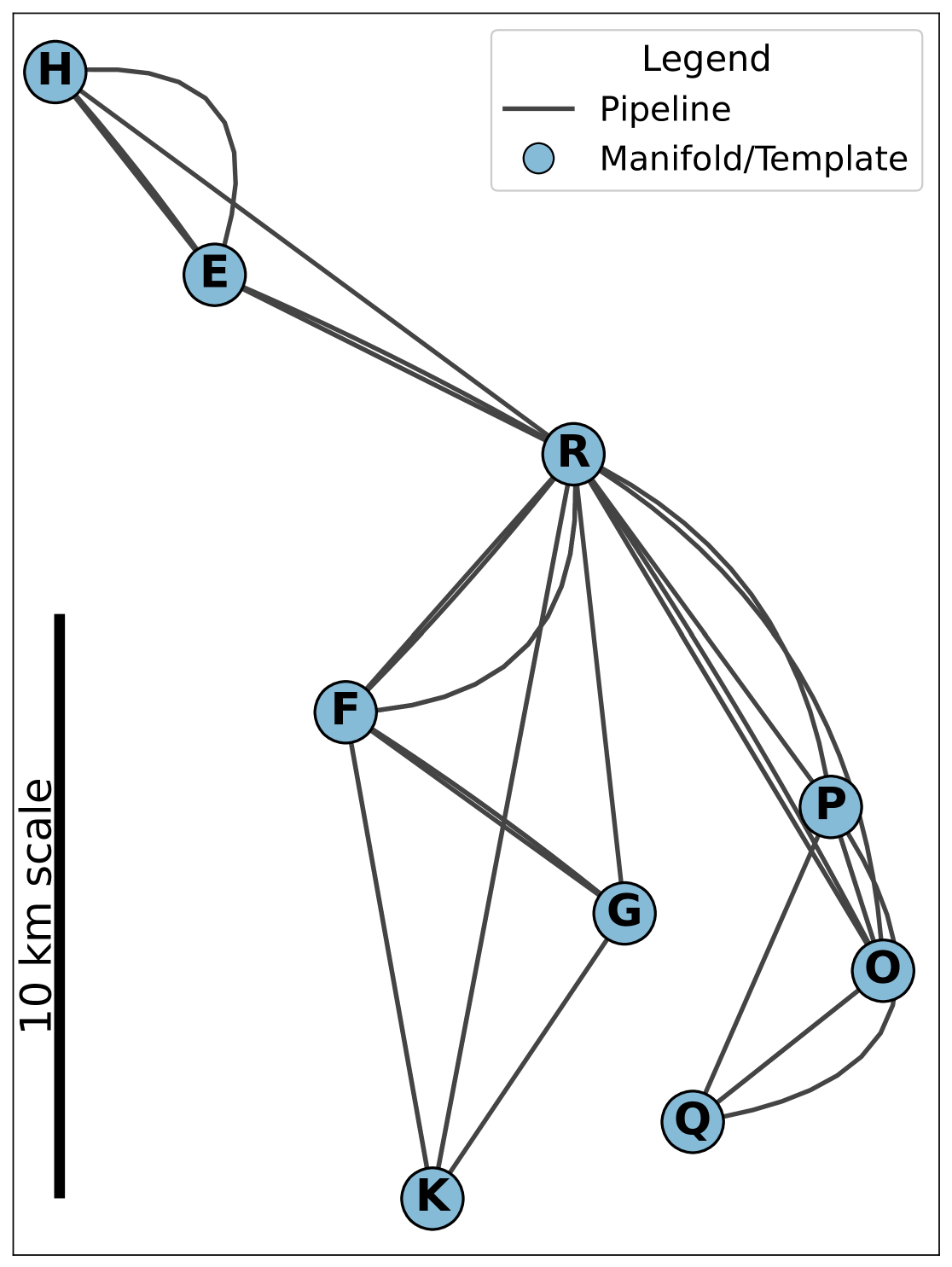}
 \caption{Map layout inspired by the Johan Sverdrup oil field \cite{norskpetroJS}, showcasing 9 structures (templates and/or manifolds) connected by 33 pipelines of varying lengths. The layout has been adapted to ensure heterogeneity and does not directly reflect proprietary or operational data for safety and confidentiality reasons.}
 \label{fig:map}
\end{figure}

This study considers a subsea oil and gas field network of pipelines and a set of candidate SDP locations (see Fig. \ref{fig:map}). Anomalies, such as leaks, dropped objects, and built-in design weaknesses, can occur anywhere along each pipeline. If a leak is detected, only the affected pipeline is identified, not the exact location of the leak. Therefore, when an alarm activates, an RAUV must:

\begin{enumerate}
  \item Depart from its home SDP,
  \item Travel to one of the two ends of the affected pipeline,
  \item Inspect along the pipeline until the anomaly is located, and
  \item Return to an SDP that minimizes the total response time.
\end{enumerate}

We define the response time as the duration between the instant an alert is received and the instant the RAUV uploads the inspection data at an SDP. The positions of the SDPs and RAUVs have a significant impact on response time. Therefore, we propose a two-phase model that minimizes both maximum and average response time, with maximum response time given higher priority. Specifically, the model determines an optimal placement of SDPs and RAUVs that minimizes the average response time (Phase 2) among those placements that minimize the maximum response time (Phase 1). For the maximum response time, we assume the anomaly (interchangeably referred to as leak) occurs at the worst-case position. For the average response time, we assume the leaks are equally likely to occur in any section of each pipeline.

The following assumptions define the operational scope of the proposed model:
\begin{enumerate}
  \item \textit{Unknown leak location.} Only the pipeline ID is known. The RAUV inspects sequentially along the pipeline until the leak is found.
  
  \item \textit{Uniform likelihood.} Anomalies are equally likely at every pipeline location. Non-uniform priors can be integrated by allocating segment-specific weights in the average response-time target, without modifying the structure of the MILP.

  \item \textit{Steady-state readiness.}  Leak events are assumed to be rare. Therefore, RAUVs are fully charged, stationed at their home SDPs, and ready for immediate deployment.

  \item \textit{Single-vehicle inspection.} Each inspection mission is executed by a single RAUV selected by the optimization model. The vehicle enters from one of the pipeline endpoints and performs a sequential inspection. If a leak is detected, the RAUV captures visual data and travels to an SDP selected by the model to upload the data. At the end of the mission, the RAUV returns to its initial SDP to be ready for another alert. Cooperative multi-vehicle inspections are not considered in the model.

  \item \textit{Station-vehicle parity.} At least as many SDPs are opened as RAUVs are deployed, with each RAUV assigned a distinct home SDP.

  \item \textit{Limited SDP locations.} In this case study, candidate SDP locations are limited to existing subsea manifolds or templates that provide the necessary power and data infrastructure. This restriction reflects the characteristics of the case study rather than a limitation of the MILP model, which can accommodate candidate SDP locations at arbitrary positions.
 
\end{enumerate}	

Solving the SDP placement problem and planning the RAUV deployment are both offline decisions. The infrastructure layout is typically fixed in a mature field, and there is no need to dynamically reallocate the RAUVs if alarms are rare events. In this context, the MILP solution time is not critical. However, it remains important to evaluate how the solution time scales as the optimization problem size increases, since the model could be used for rapid scenario analyses and to assess larger operational areas, for instance involving multiple fields of the Norwegian Continental Shelf.

\section{SOLUTION METHOD}
\label{solution_method}

A set of subsea pipelines 
$\mathcal{P}$ requires incident inspection whenever an anomaly is detected. There are $\numCandidateStations$ candidate SDP locations, from which exactly $\numChosenStations$ sites are chosen (with $\numChosenStations \le \numCandidateStations$) to host RAUVs. The RAUV fleet consists of $\numRobots$ vehicles, where $\numRobots\le\numChosenStations$.

Readings of the pressure of the pipelines can identify a potential leak in a particular pipeline, but cannot locate the potential leak. If a potential leak is detected in the pipeline $p$, one of the RAUVs residing at the stations travels to one end of the pipeline, inspects the pipeline until it finds the leak, and collects visuals, and travels from the leak to the station to upload the visuals for the operations support team.

The goal is to select a set of SDPs and assign each RAUV to one of them as its home base, in order to minimize the time required to inspect any pipeline anomaly. From here on, the terms station and SDP will be utilized interchangeably, as will the terms robot and RAUV. Each inspection mission consists of three sequential phases:
\begin{enumerate}
 \item \textbf{Approach phase:}
  The RAUV departs its home SDP and travels to one of the two endpoints of the target pipeline at a transit speed.
 \item \textbf{Inspection phase:} 
  From that endpoint the RAUV moves along pipeline toward the anomaly, at an inspection speed.
 \item \textbf{Return phase:} 
  After reaching the anomaly, the RAUV returns to an open SDP that minimizes the total response time, again at transit speed. 
\end{enumerate}

\subsection{Index Sets}
The following information describes the sets utilized in the MILP. 
\begin{itemize}
 \item $\mathcal{P}=\{1,\dots,\numPipelines\}$: pipelines, indexed by $p$.
 \item $\mathcal{Q}=\{1,\dots,\numCandidateStations\}$: candidate SDP locations, indexed by $q$, see Table \ref{tab:station_map} for the mapping.
\end{itemize}

\begin{table}[ht]
\begin{center}
\caption{Mapping from station index $q$ to its letter label.}
\label{tab:station_map}
\begin{tabular}{| c | c |}
\hline
\textbf{Index $q$} & \textbf{Label} \\ 
\hline
1 & H \\ 
\hline
2 & E \\ 
\hline
3 & R \\ 
\hline
4 & F \\ 
\hline
5 & G \\ 
\hline
6 & K \\ 
\hline
7 & O \\ 
\hline
8 & P \\ 
\hline
9 & Q \\ 
\hline
\end{tabular}
\end{center}
\end{table}

\subsection*{Parameters}
\begin{itemize}
 \item $\numChosenStations$: number of SDPs to select, $\numChosenStations\le\numCandidateStations$. 
 \item $\numRobots$: number of RAUVs to deploy, $\numRobots\le\numChosenStations$.
 \item $\transitSpeed$: transit speed of each RAUV (m/s).
 \item $\inspectSpeed$: inspection speed along the pipeline (m/s).
 \item $\numSamples{p}$: number of uniform sample sections on pipeline $p$.
 \item $\anomalyCoord$: Cartesian coordinates of the $i$th potential anomaly sample (i.e., pipeline section) on pipeline $p$.
  \item $\entryCoord{p,\delta}$: coordinates of the two endpoints of pipeline $p$. $\delta=0$ indicates the endpoint in section 1 and $\delta=1$ indicates the endpoint in section $N_p$, that is $\entryCoord{p,0} = \mathbf{a}_{p, 1}$ and $\entryCoord{p,1} = \mathbf{a}_{p, N_p}$.
 \item $\stationCoord{j}$: Cartesian coordinates of a candidate SDP $j$.
 \item $\epsParam$: allowable tolerance on maximum time in Phase 2. It is used to improve solver convergence.
\end{itemize}

\subsection{Geometric Preprocessing}
The proposed approach divides each pipeline $p$ of length $L_p$ into $N_p$ segments of equal length $\Delta s$ except for the last section that might be shorter, with $L_p \lesssim N_p \, \Delta s$. The time to inspect a section is:

\begin{equation}
 \Delta t 
 = \frac{\Delta s}{\inspectSpeed}.
\end{equation}

The following definitions then apply:

\begin{itemize}
 \item \textbf{Approach time:}
  \begin{equation}  
  t^{\rm app}_{p,\delta,j} 
  = \frac{\|\entryCoord{p,\delta} - \stationCoord{j}\|}{\transitSpeed}.
  \end{equation}
 \item \textbf{Inspection time (from one of the pipeline's endpoints to an anomaly $i$):}
  \begin{equation} 
  t^{\rm ins}_{p,\delta,i}  = 
    \begin{cases}
      i\,\Delta s / \inspectSpeed, & \text{if } \delta=0 \\
      (N_p-i+1)\,\Delta s / \inspectSpeed, & \text{if } \delta=1
    \end{cases}.
  \end{equation}

 \item \textbf{Return time:}
  \begin{equation}
  t^{\rm ret}_{p,i,k} 
  = \frac{\|\anomalyCoord - \stationCoord{k}\|}{\transitSpeed}.
  \end{equation}
\end{itemize}

\subsection{MILP Formulation}
The MILP formulation includes decision variables and constraints that limit the outcomes.  
\paragraph{Decision variables}
\begin{itemize}
 \item $s_j\in\{0,1\}$ for $j\in\mathcal Q$: 1 if SDP $j$ is available (in other words, an SDP is positioned there).
 \item $r_j\in\{0,1\}$ for $j\in\mathcal Q$: 1 if an RAUV is based at station $j$.
 \item $x_{p,\delta,j}\in\{0,1\}$, $p\in\mathcal P$, $\delta\in\{0,1\}$, $j\in\mathcal Q$: 1 if pipeline $p$ is entered at endpoint $\delta$ by RAUV from $j$, with $\delta=0$ if the endpoint is  section 1 of the pipeline and $\delta=1$ if the endpoint is section $N_p$ of the pipeline.
 \item $y_{p,i,k}\in\{0,1\}$, $p\in\mathcal P$, $i=1,\dots,N_p$, $k\in\mathcal Q$: 1 if after inspecting segment $i$ on $p$, RAUV returns to $k$.
\end{itemize}

\paragraph{Station \& RAUV budgets}
\begin{equation}
\label{eq:budget}
 \sum_{j\in\mathcal Q} s_j = S,
 \quad
 \sum_{j\in\mathcal Q} r_j = R,
 \quad
 r_j \le s_j,
 \quad\forall j\in \mathcal Q.
\end{equation}

Equation \eqref{eq:budget} enforces three constraints. The first, $\sum_{j \in \mathcal{Q}} s_j = S$, ensures that exactly $S$ of the $Q$ candidate SDPs are selected. The second, $\sum_{j \in \mathcal{Q}} r_j = R$, specifies that exactly $R$ RAUVs are deployed across those stations. Finally, the inequality $r_j \le s_j$ guarantees that an RAUV can only be based at station $j$ if that station is made available.

\paragraph{Entry assignment}
\begin{align}
 &\sum_{\delta=0}^1\sum_{j\in\mathcal Q} x_{p,\delta,j} = 1,
  &&\forall p\in\mathcal P, 
  \label{eq:entry1}\\
 &x_{p,\delta,j} \le r_j,
  &&\forall p\in\mathcal P,\;\delta\in\{0,1\},\;j\in\mathcal Q.
  \label{eq:entry2}
\end{align}

Constraint (\ref{eq:entry1}) forces each pipeline $p$ to choose exactly one endpoint $\delta$ and departure from station $j$. The inequality $x_{p,\delta,j}\le r_j$ ensures that entry can only occur from stations where an RAUV is actually based.

\paragraph{Return assignment}
\begin{align}
 &\sum_{k\in\mathcal Q} y_{p,i,k} = 1,
  &&\forall p\in\mathcal P,\;i=1,\dots,N_p,
  \label{eq:return1}\\
 &y_{p,i,k} \le s_k,
  &&\forall p\in\mathcal P,\;i=1,\dots,N_p,\;k\in\mathcal Q.
  \label{eq:return2}
\end{align}

Here, \eqref{eq:return1} ensures that for each pipeline $p$ and anomaly $i$, exactly one return station $k$ is selected. Constraint \eqref{eq:return2} enforces that returns can only be assigned to stations that have been made available ($s_k=1$). 

\paragraph{Maximum (Phase 1) cost}
For every pipeline $p$ and possible leak position $i$:

\begin{equation}
\begin{aligned}
&\sum_{\delta=0}^1 \sum_{j \in \mathcal{Q}}
\big( t^{\mathrm{app}}_{p,\delta,j} + t^{\mathrm{ins}}_{p,\delta,i} \big)
\, x_{p,\delta,j} \\
&\quad+\; \sum_{k \in \mathcal{Q}} t^{\mathrm{ret}}_{p,i,k} \, y_{p,i,k}
\;\le\; T_{\max}, \quad 
\forall p \in \mathcal{P},\ \forall i \in \mathcal{N}_p.
\end{aligned}
\label{eq:phase1}
\end{equation}

This constraint bounds the total mission time by $T_{\max}$, enforcing a common maximum limit across all possible leak locations.

\paragraph{Average (Phase 2) cost}

In Phase 2, we minimize the average inspection response time across all possible leak locations:

\begin{align}
\label{eq:phase2_avg}
\frac{1}{\sum_{p \in \mathcal{P}} N_p}
\sum_{p \in \mathcal{P}} \sum_{i=1}^{N_p} \Bigg(
& \sum_{\delta=0}^1 \sum_{j \in \mathcal{Q}}
\big( t^{\mathrm{app}}_{p,\delta,j} + t^{\mathrm{ins}}_{p,\delta,i} \big)
\, x_{p,\delta,j} \notag \\
& + \sum_{k \in \mathcal{Q}}
t^{\mathrm{ret}}_{p,i,k} \cdot y_{p,i,k}
\Bigg)
\;\le\;
T_{\text{avg}},
\end{align}

with the linking constraint:
\begin{equation}
\label{eq:phase2_link}
T_{\max} \le (1+\epsilon) T_{\max}^\star,
\end{equation}
where $T_{\max}^\star$ is the Phase 1 optimum and $\epsilon$ a small tolerance (e.g., $10^{-3}$).
This ensures $T_{\max}$ in Phase 2 stays close to its best case value from Phase 1, preventing large maximum case delays while improving the average. The two‐stage design thus balances both objectives.

The MILP formulation presented above is implemented and solved using the \texttt{PySCIPOpt} optimization framework~\cite{MaherMiltenbergerPedrosoRehfeldtSchwarzSerrano2016}, which provides a Python interface to SCIP, a solver renowned for efficiently handling large scale integer and combinatorial optimization problems~\cite{Achterberg2009}. 

\section{RESULTS}
\label{results}

The results section details the performance of the MILP under different combinations of RAUVs ($\numRobots$) and chosen SDPs ($\numChosenStations$), resulting in 45 scenarios (combinations of $\numRobots=1,2,\ldots,9$ and $\numChosenStations=\numRobots,\numRobots+1,\ldots,9$). The two key performance indicators are the maximum and average response times on a map layout inspired by the Johan Sverdrup oil and gas field consisting of 9 candidate SDP positions and 33 pipelines. The following figures are used to examine performance, SDP and RAUV base usage patterns, and the trade-off between performance and deployment cost. 

\subsection{Impact of Infrastructure vs Fleet}

\begin{figure}[ht]
\centering
\includegraphics[width=\columnwidth]{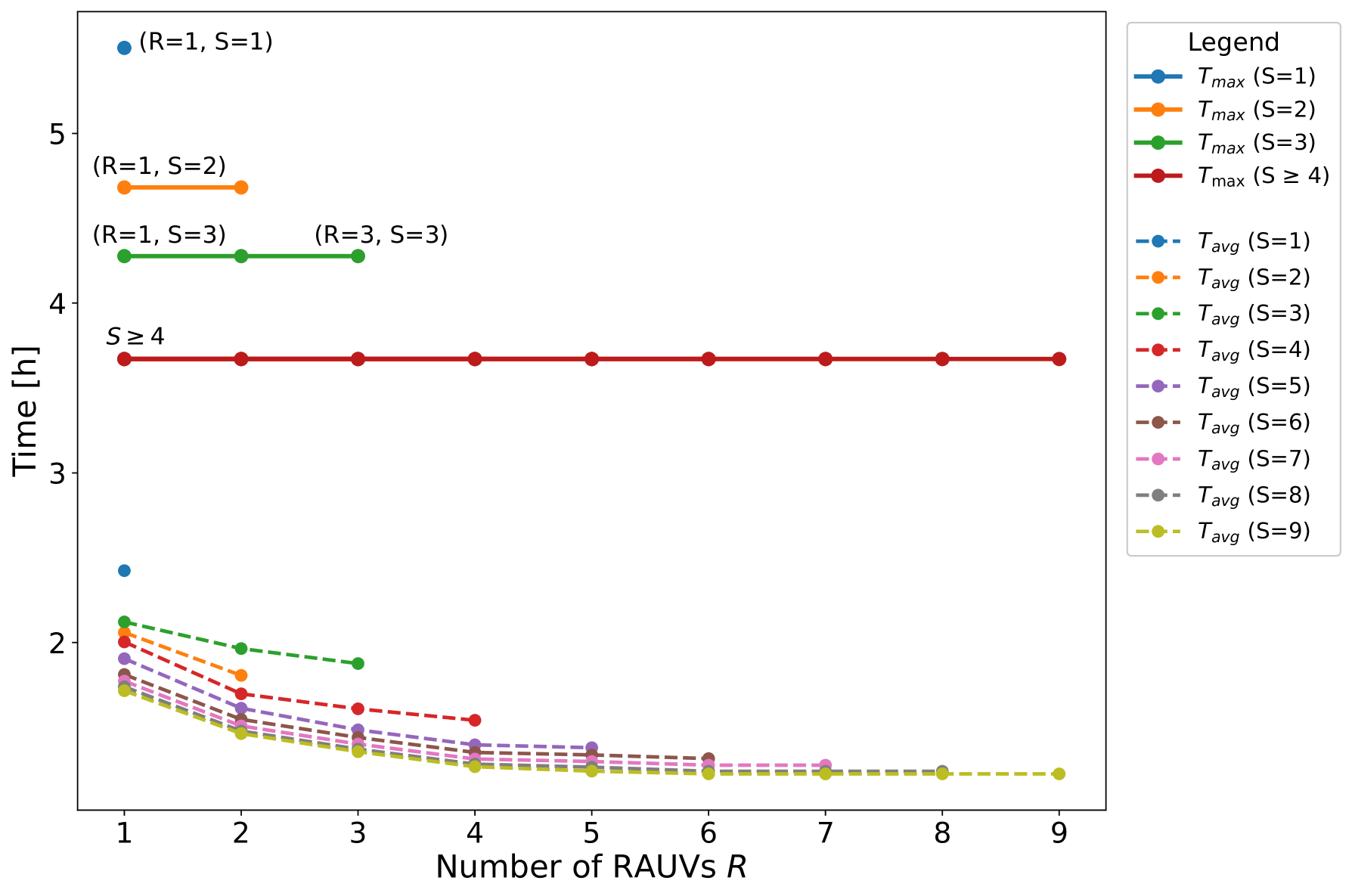}
\caption{Maximum ($\tMax$, solid lines) and average ($\tAvg$, dashed lines) response times as a function of the number of RAUVs $\numRobots$, for different numbers of SDPs $\numChosenStations$. Each color pair (solid/dashed) corresponds to a fixed $\numChosenStations \in \{1,\dots,9\}$. For $S \geq 4$, the maximum response time converges to approximately 3.6~h across all fleet sizes.}
\label{fig:x_robots}
\end{figure}

Fig.~\ref{fig:x_robots} plots max response time $\tMax$ (solid line) and average response time $\tAvg$ (dashed line) against fleet size $\numRobots$, for a given budget on the number of stations  $\numChosenStations$. The results reveal the difference between the worst-case versus the average case response time.

The max response time $\tMax$ is dominated by SDP placement rather than fleet size. Increasing the fleet size has little impact on $\tMax$, whereas increasing the number of SDPs can drastically reduce the worst-case mission time by optimizing spatial coverage. When the spatial coverage is sufficient (four SDPs in the figure), $\tMax$ levels off at around 3.6 hours, indicating that the worst-case pipeline can be served from a well-positioned station.

The average case, on the other hand, continues to decrease as both SDPs and RAUVs increase in number. More stations reduce the travel time, and the increased fleet size provides more flexibility in the assignments, further reducing the average time.

The convergence behavior of $\tAvg$ aligns with the underlying geometry of the field. A theoretical lower bound of 1.22~h was computed by assuming that each pipeline is entered from its optimal endpoint and returned to the closest SDP, under idealized conditions with no constraints on fleet size. The MILP solutions closely approach this limit. For example, a configuration with $R=5$, $S=9$ achieves $\tAvg = 1.24$~h, which is only $\approx$1.6\% above the theoretical minimum under idealized entry and return conditions. Moreover, 5\% convergence is already reached by $R=4$, $S=8$ with $\tAvg = 1.28$~h.

This behavior reflects spatial saturation: by $S = 5$, 99.6\% of leak points lie within 2~km of at least one SDP, and 99\% have two SDPs within 3~km. At this point, both entry and return travel times approach their minimum possible values, and increasing the number of stations beyond this provides only marginal reductions in response time.

In summary, near-optimal maximum and average response times can be achieved with only a few well-placed RAUV–SDP pairs. Configurations with $R=4-5$ and $S=6-8$ already yield response times within 1--5\% of the minimum, underscoring the efficiency of the optimization approach and the saturation of the spatial environment.

\subsection{Station Usage Analysis: Selected vs. Home Roles}

To understand the geographic patterns in the MILP solutions, it is first important to examine how often each candidate SDP $q\in\mathcal{Q}$ is included in the chosen set $\mathcal{S}$ of size $\numChosenStations$. 

\begin{figure}[ht]
\centering
\includegraphics[width=0.8\columnwidth]{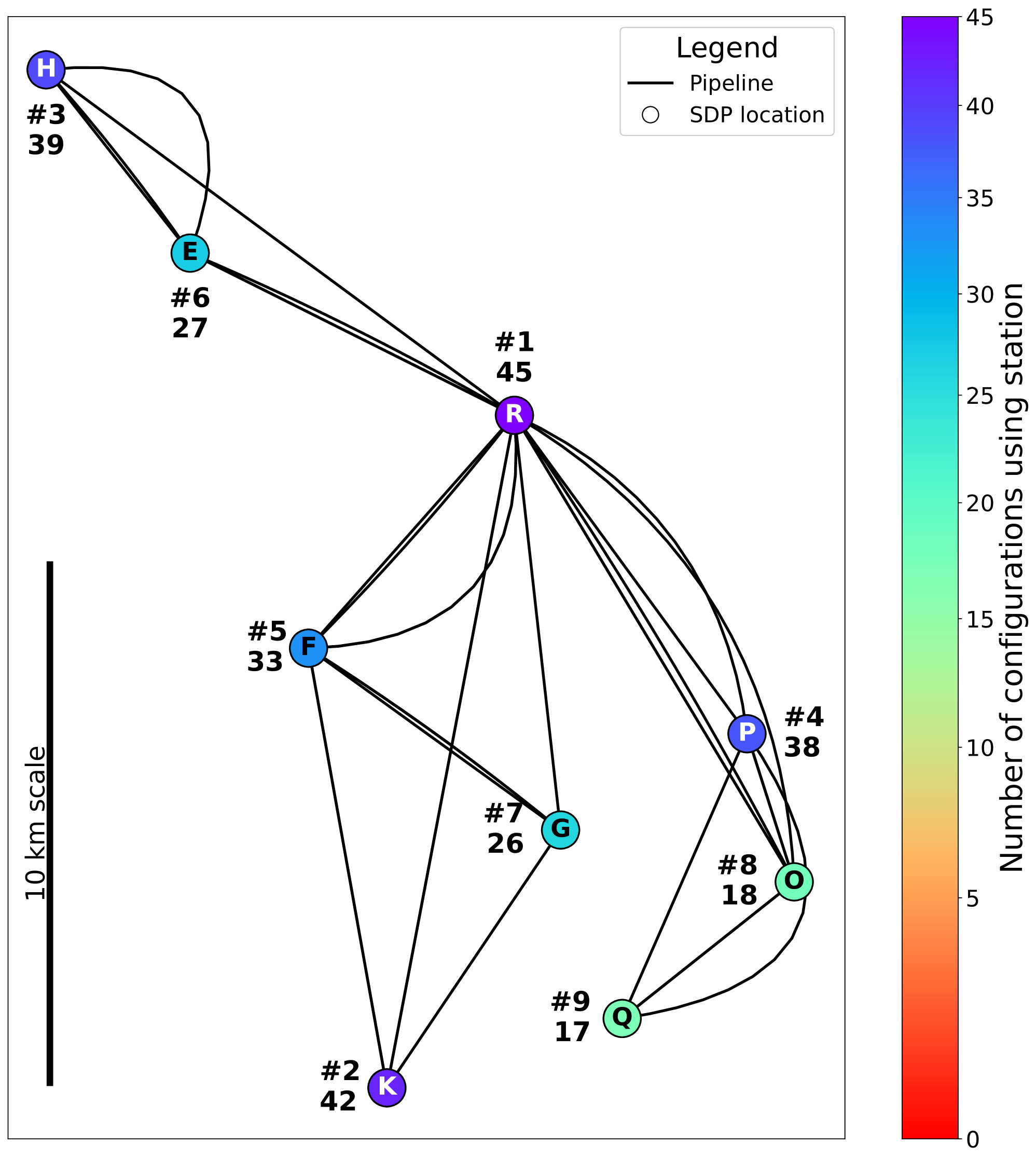}
\caption{Frequency with which each SDP is included in the selected set $\mathcal{S}$ across all $(\numRobots,\numChosenStations)$ scenarios (45 total). Station rank (\#1 = most frequent) is shown next to each SDP location.}
\label{fig:heatmap_stations}
\end{figure}

Fig.~\ref{fig:heatmap_stations} colors each node by the number of $(\numRobots,\numChosenStations)$ configurations (out of 45 total) for which $q\in\mathcal{S}$. Station R is used in all 45 configurations, followed closely by station K (42 selections) and station H (39). These three locations form the core of the optimal deployment strategy across nearly all tested fleet and station budgets.

Their consistent inclusion reflects their topological roles. Station R lies at the center of the network and connects efficiently to the majority of pipelines, particularly those to the north and east. Station K, located in the southern region, provides necessary access to a geographically distant cluster of pipelines, while station H serves as a northern anchor connecting pipelines near station E.

Mid-frequency stations include P (38), F (33), E (27), and G (26). These stations tend to complement the core sites when larger station budgets allow for more precise access to pipelines. For instance, P and F provide access to central and eastern pipelines, while E and G cover the north-west and central-south, respectively. Despite being moderately used, their presence can improve average response time and allow for more flexibility, particularly in scenarios with $S \ge 6$.

Low frequency stations include O (18) and Q (17), which appear in fewer than half the configurations. Their sporadic inclusion suggests that these stations are only selected when higher station budgets afford redundancy.

Importantly, spatial proximity does not always imply similar usage rates. Although stations P and O are geographically close, P is used in 38 configurations, while O appears in just 18. Similarly, E is selected in 27 configurations while its nearby neighbor H appears in 39. These discrepancies underscore that pipeline connectivity and strategic positioning dominate over raw geometric placement in determining station value.

\begin{figure}[ht]
\centering
\includegraphics[width=0.8\columnwidth]{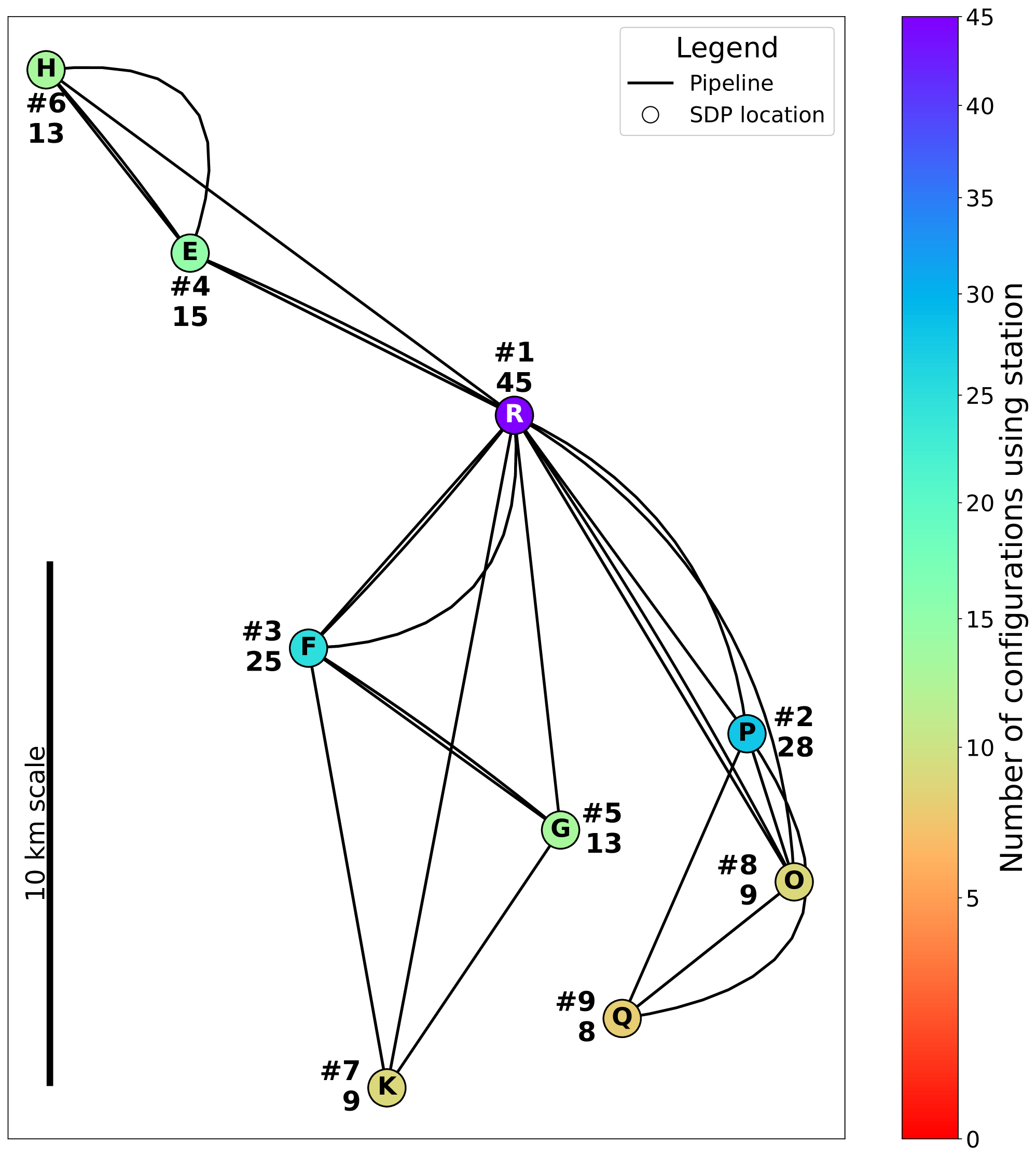}
\caption{Frequency with which each SDP serves as a home station for an RAUV across all $(\numRobots,\numChosenStations)$ scenarios. Station rank (\#1 = most frequent) is shown next to each SDP location.}
\label{fig:heatmap_home_stations}
\end{figure}

Fig.~\ref{fig:heatmap_home_stations} shows the number of configurations in which each station $q$ serves as a \textit{home} base for at least one RAUV. Here, station R again dominates, serving as a home location in all 45 scenarios. The next most frequent home stations are P (28 configurations) and F (25). 

Compared to station selection frequencies, the divergence is notable. Station K, for example, appears in 42 of 45 SDP selections but serves as a home base in only 9 cases. This suggests that while K is often useful to facilitate return paths or reduce average response time, it is suboptimal as a starting location due to its relatively outer position in the field. A similar pattern holds for stations H and Q, which are used as home locations in only 13 and 8 scenarios respectively.

In contrast, stations like F and E are preferred as launch points due to their more central or strategically balanced locations relative to the pipeline network. Interestingly, some spatially adjacent stations show strong divergence in home usage: P is a home base in 28 configurations, whereas its neighbor O is chosen in only 9. 

These findings emphasize that home station selection is more selective than general SDP usage. Only a small subset of SDPs regularly serve as reliable RAUV launch sites, while others contribute mainly to minimizing return time or supporting specific configurations. Prioritizing basing at stations with broad pipeline access enables the system to maintain strong performance with minimal redundancy in vehicle placement.

\subsection{Cost vs. Response Time Trade-off}

\begin{figure}[ht]
\centering
\includegraphics[width=0.48\textwidth]{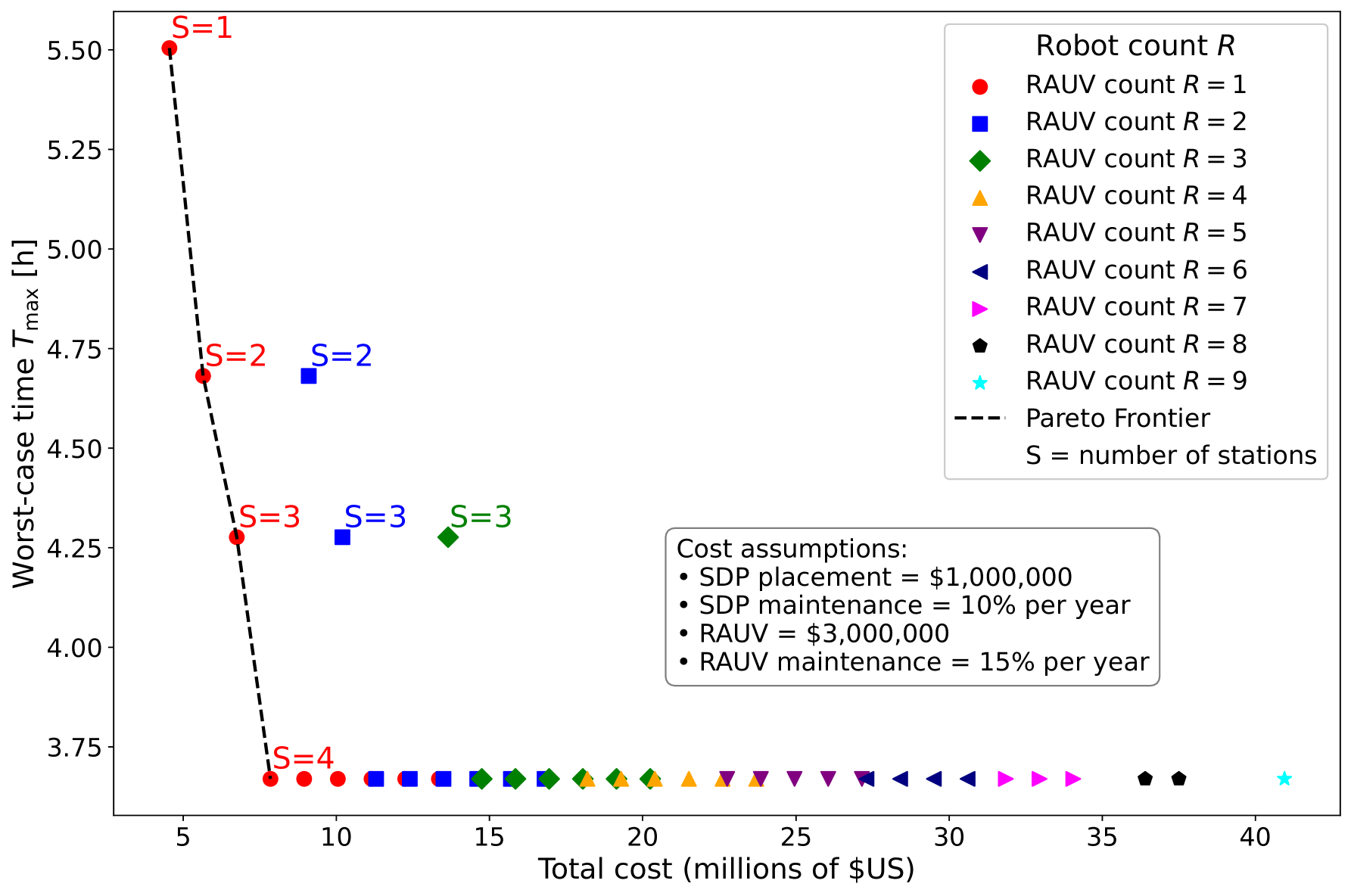}
\caption{Trade‐off between total deployment cost (in millions of USD) and maximum response time $\tMax$. The dashed black line is the Pareto frontier.}
\label{fig:max_cost}
\end{figure}

Fig.~\ref{fig:max_cost} plots total investment in USD (horizontal axis) against the maximum response time $\tMax$ in hours (vertical axis) for each feasible $(\numRobots,\numChosenStations)$ configuration. The total cost is computed as

\begin{equation}
\label{eq:cost_rewritten}
\begin{split}
\text{Cost}
 &= S\,C_{\mathrm{SDP}}\bigl(1 + m_{\mathrm{SDP}}\bigr) \\
 &\quad
  +\,R\,C_{\mathrm{RAUV}}\bigl(1 + m_{\mathrm{RAUV}}\bigr),
\end{split}
\end{equation}

where $C_{\mathrm{SDP}}$ and $C_{\mathrm{RAUV}}$ are the unit placement costs of an SDP and an RAUV, respectively, and $m_{\mathrm{SDP}}, m_{\mathrm{RAUV}}\in[0,1]$ are their first year maintenance fractions. The cost assumptions are listed in the figure's legend. Each point in Fig.~\ref{fig:max_cost} represents a unique $(R,S)$ configuration, with the horizontal position indicating total cost and the vertical position the corresponding $\tMax$ in hours. Labels of the form “$S=\dots$” indicate the number of SDPs used in that configuration.

A clear Pareto frontier emerges (dashed black line), delineating the most cost-effective configurations. From the leftmost point at $(R=1, S=1)$, where $\tMax \approx 5.5$~h and the cost is \$4.55 million, increasing the number of stations to two reduces $\tMax$ to approximately 4.6~h at an additional cost of \$1.1 million. Adding a third station ($S=3$) further reduces the maximum response time to 4.2~h for a total of \$6.75 million. The final configuration of $R=1, S=4$, marking the end of the frontier drops the maximum response time to 3.6~h at a cost of \$7.8 million. Beyond this configuration, adding RAUVs or SDPs yields no additional gains: $\tMax$ remains constant at 3.6~h. All configurations above the frontier are not cost-effective as they incur higher costs without improving maximum performance.

\begin{figure}[ht]
\centering
\includegraphics[width=0.48\textwidth]{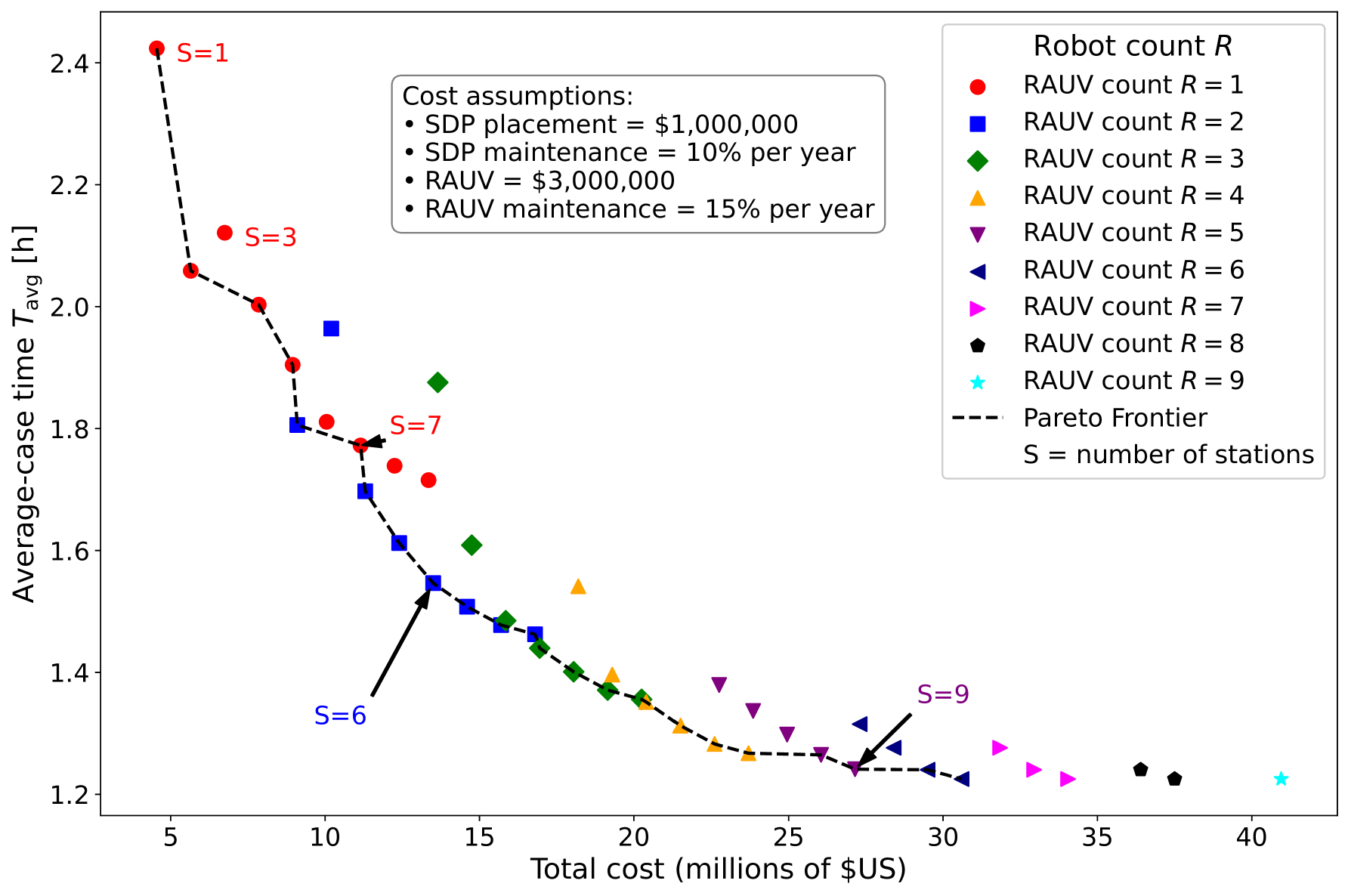}
\caption{Trade‐off between total deployment cost (in millions of USD) and average case response time $\tAvg$. The dashed black line is the Pareto frontier.}
\label{fig:avg_cost}
\end{figure}

Fig.~\ref{fig:avg_cost} presents a similar trade-off for the average case response time $\tAvg$. The shape of the frontier follows a smooth elbow pattern. The configuration $(R=1, S=1)$ yields the highest average response time on the frontier, $\tAvg$ = 2.4~h at a cost of \$4.55 million. Adding a second station reduces this to 2.06~h at \$5.65 million. The best 1 RAUV setup is $(R=1, S=7)$, achieving 1.77~h at a cost of \$11.1 million. The figure shows that incorporating additional RAUVs leads to more cost efficient improvements.

A nuanced finding emerges when considering a third station. The configuration $(R=1,S=3)$ results in a slightly increased average response time of $\tAvg$ = 2.12~h, even though the budget increases to \$6.75 million. This indicates a trade-off where the model, in seeking a better overall solution (specifically for max-case performance), selects a set of stations that, while more robust, are less efficient for the average case.

The average time continues to decrease with further investments. The configuration $(R=2,S=6)$ achieves $\tAvg$ $\approx$ 1.55~h at a cost of \$13.50 million, marking a key point where the cost-to-performance ratio becomes less favorable. The lowest average time on the frontier is $\tAvg$ $\approx$ 1.23 h, achieved by the configuration $(R=6,S=9)$ at a cost of \$30.60 million.

The analysis reveals a significant difference in the performance and cost trade-off between the two objectives. While minimizing the maximum response time saturates at a lower budget, minimizing the average time below a certain threshold (e.g., 1.5 h) requires a significantly higher investment. For example, a budget of approximately \$15 million for the $(R = 2, S = 7)$ configuration balances deployment cost and inspection performance with a $\tAvg$ $\approx$ 1.51 h.

\subsection{Scalability Analysis: MILP Solve Times}

\begin{figure}[ht]
 \centering
 \includegraphics[width=0.48\textwidth]{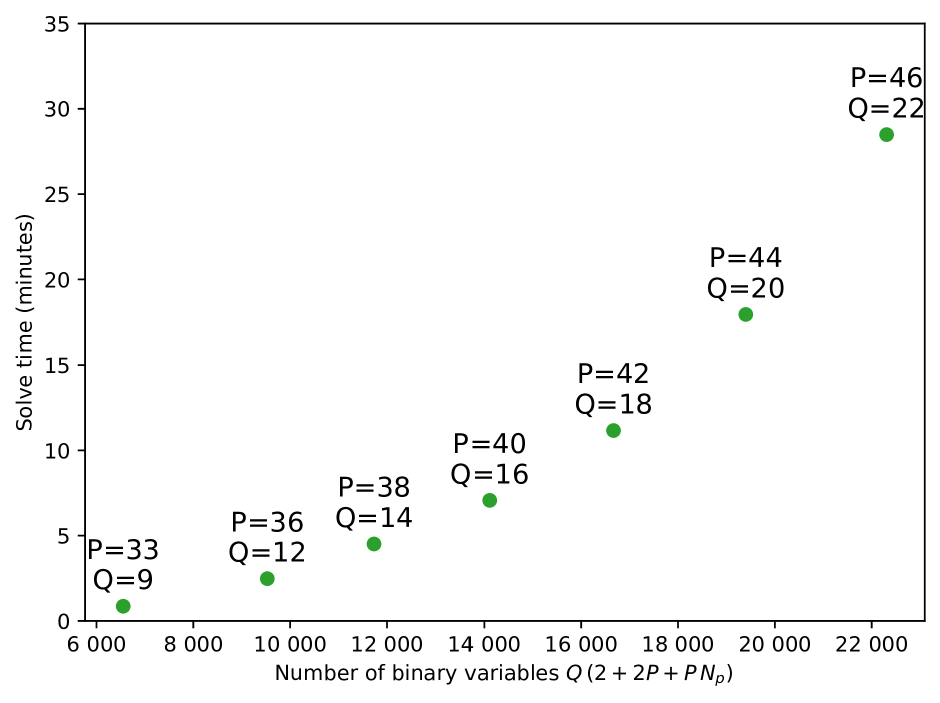}
 \caption{Average two-phase MILP solve time (in minutes) as a function of the total number of binary decision variables. Each point is labeled with the corresponding number of pipelines $P$ and candidate stations $Q$.}
 \label{fig:scaling_curve}
\end{figure}

The number of binary variables in the MILP grows with the number of candidate stations $Q$, pipelines $P$, and leak samples $N_p$, scaling as
\begin{equation}
2Q + 2PQ + Q\sum_{p=1}^{P} N_p,
\end{equation}
which is $\mathcal{O}(PQ N_{\max})$ in the maximum case. As the map expands, so does solve time: Fig.~\ref{fig:scaling_curve} shows that doubling variable count increases the runtime sixfold.

\begin{figure}[ht]
 \centering
 \includegraphics[width=0.48\textwidth]{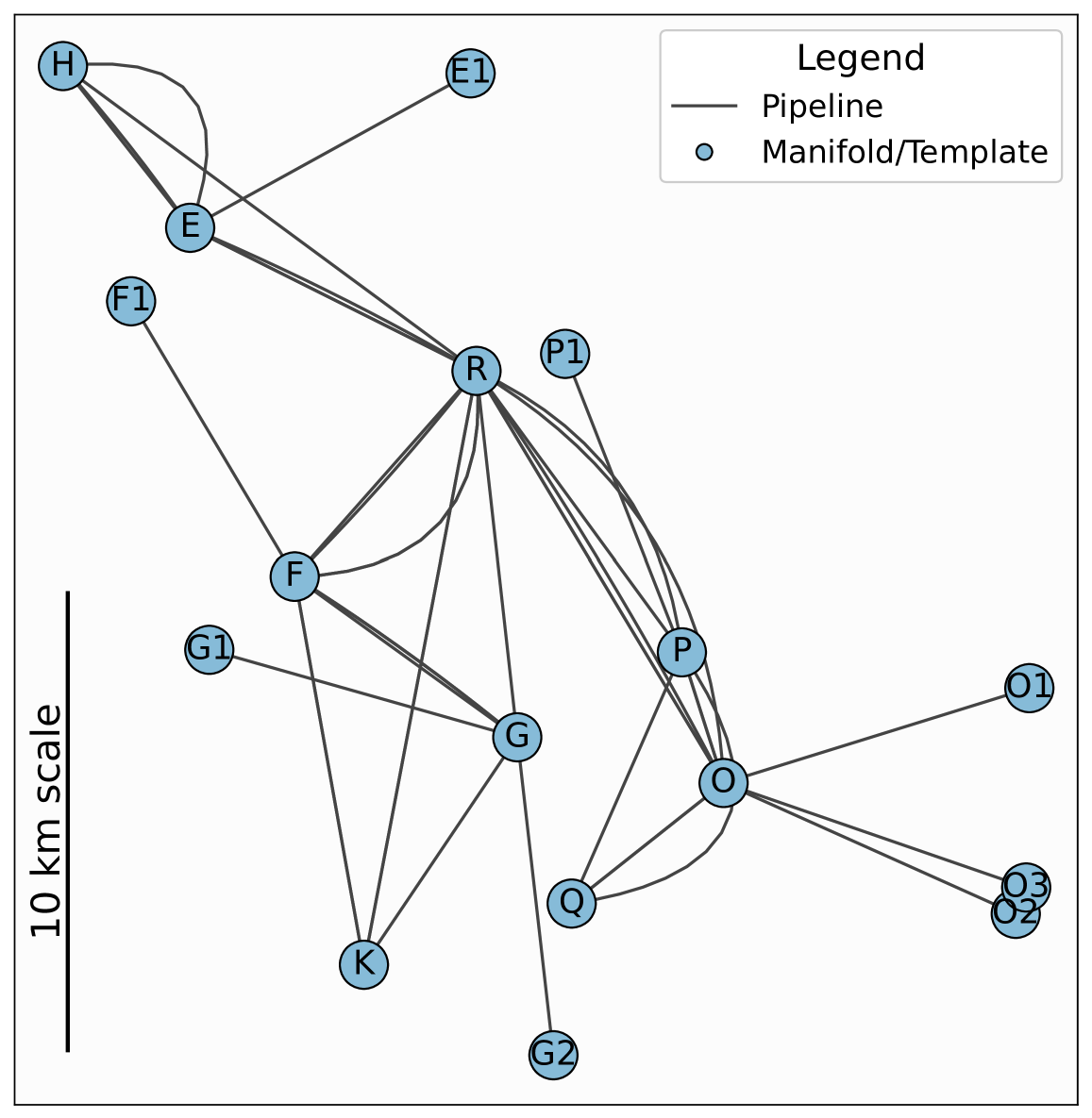}
 \caption{An expanded Johan Sverdrup map with 8 additional candidate SDP locations and pipelines set at random, resulting in a map with 17 candidate SDP locations and 41 pipelines.}
 \label{fig:fake_map}
\end{figure}

Experiments on synthetic maps (e.g., Fig. ~\ref{fig:fake_map}) show that while small maps ($Q=9$, $P=33$) solve all 45 scenarios in just over a minute, larger maps ($Q=22$, $P=46$) take over 87 minutes per scenario on average. Table~\ref{tab:scalability} summarizes solve-time statistics for each map size.

A consistent pattern is observed: configurations with $R=1$ and low $S$ are hardest to solve. For example, at $(Q=22, P=46)$, the $(R=1, S=2)$ case takes 87 minutes, while all other configurations remain under 11 minutes. These low resource cases likely challenge the solver due to tighter constraints and large combinatorial search space. 

\begin{table}[htbp]
 \centering
 \caption{Solve‐time summary for $N_p=20$ on a 32-core workstation.}
 \label{tab:scalability}
 \resizebox{\columnwidth}{!}{%
  \begin{tabular}{|c|c|c|c|c|c|}
   \hline
   $\boldsymbol{Q}$ 
   & $\boldsymbol{P}$ 
   & $\text{Num. binaries}$ 
   & $\text{Avg.\ solve (min)}$ 
   & $\text{Num. scenarios}$ 
   & $\text{Total (min)}$ \\ 
   \hline
   9  & 33 & 6 552 & 0.86 & 45 & 1.21 \\ \hline
   12 & 36 & 9 528 & 2.48 & 78 & 6.04 \\ \hline
   14 & 38 & 11 732 & 4.51 & 105 & 14.81 \\ \hline
   16 & 40 & 14 112 & 7.07 & 136 & 30.05 \\ \hline
   18 & 42 & 16 668 & 11.16 & 171 & 59.65 \\ \hline
   20 & 44 & 19 400 & 17.94 & 210 & 117.75 \\ \hline
   22 & 46 & 22 308 & 28.49 & 253 & 225.23 \\ \hline
  \end{tabular}
 }
\end{table}

\section{DISCUSSION}
\label{discussion}

\subsection{System Trade‐Offs}

\begin{figure}[ht]
\centering
\includegraphics[width=0.8\columnwidth]{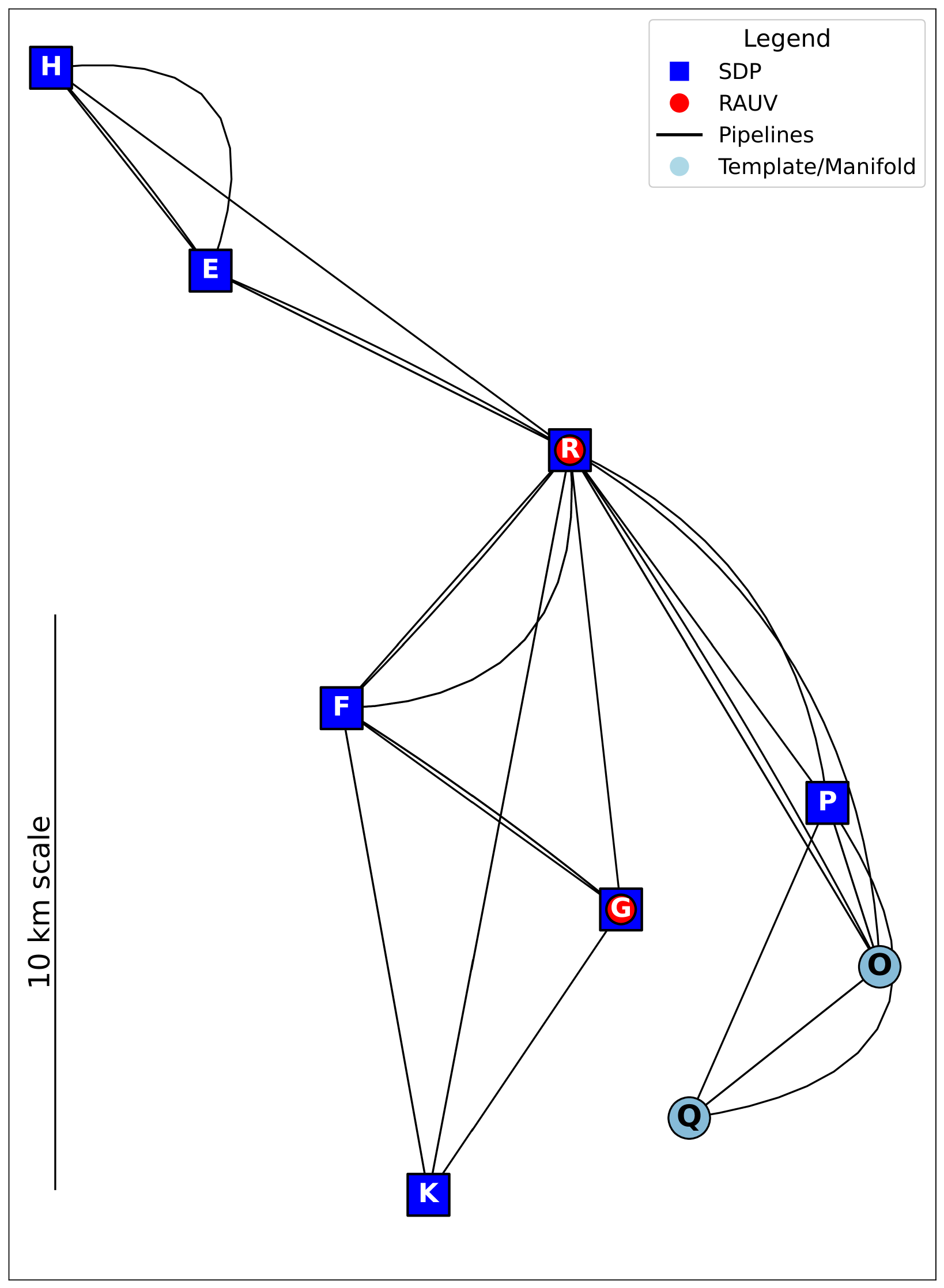}
\caption{The map depicts an example of optimal placement for two RAUVs and seven SDPs, with the RAUVs situated at R and G.}
\label{fig:best_example_map}
\end{figure}

The trade-off analysis reveals that a configuration with 2 RAUVs and 7 SDPs offers a compelling balance between performance and cost. This setup achieves a maximum response time of $\tMax = 3.6$~h and an average case time of $\tAvg = 1.5$~h both near the minimum of the Pareto frontiers while keeping total cost (\$14.6 million) well below the point where diminishing returns become severe. As shown in Fig.~\ref{fig:best_example_map}, the recommended assets are strategically placed to cover the field efficiently with minimal redundancy.

More generally, the results show a clear pattern of diminishing marginal returns as the number of RAUVs increases, especially beyond $\numRobots \approx 5$. In Fig.~\ref{fig:x_robots}, both $\tMax$ and $\tAvg$ exhibit steep improvements when increasing from one to four RAUVs, but further additions yield only minor gains.

By contrast, increasing the number of SDPs from $S=1$ to $S=4$ even with a single RAUV produces substantial improvements. These early-stage gains dominate the cost-performance trade-off.

The Pareto frontier in Fig.~\ref{fig:avg_cost} quantifies the performance and cost trade-off. While the frontier as a whole shows a decrease in average response time as investment increases, this is not true for all individual configurations. For example, adding a third station to a single RAUV fleet counterintuitively results in a small increase in average response time, from 2.06~h to 2.12~h (3.5 min). This demonstrates that increasing assets does not guarantee improved average performance, as the optimal placement for reducing maximum response time can be less efficient for the average case.

In summary, the analysis demonstrates that a small and well-placed fleet, particularly with a few central SDPs can achieve low inspection times without excessive investment. Operators can use these findings to plan deployments that align with performance requirements, infrastructure constraints, and budget ceilings.

\subsection{Spatial Insights}

Figs.~\ref{fig:heatmap_stations} and~\ref{fig:heatmap_home_stations} show how often each candidate SDP is included in the selected set $\numChosenStations$ or assigned as a home base across all $(R,S)$ configurations. The first figure reports the number of scenarios in which each station is selected as part of the active deployment, while the second counts how frequently each location hosts an RAUV.

Stations R, K, and H are selected in the majority of configurations, reflecting their geometric advantage in reducing overall mission time. Their consistent presence suggests they serve as critical access points to large portions of the pipeline network. While stations like P and F are also frequently chosen, more peripheral locations such as Q and O appear less often, typically only when larger station budgets are allowed.

Home station assignments diverge from selection frequency. For example, K is chosen as an SDP in 42 configurations but serves as a home station in only 9, implying it plays a supporting role for return routing rather than deployment. In contrast, R, P, and F appear as frequent home bases, indicating that optimal placement can often be achieved by anchoring RAUVs at just a few well-placed hubs.

Placing SDPs is largely driven by geometric necessity for pipeline access, while placing RAUVs is much more nuanced: it requires balancing entry distance, flexibility across missions, and overall time efficiency.

This pattern underscores a key spatial insight: effective mission performance does not require a uniform distribution of RAUVs across all SDPs. Instead, strategic placement at high-impact sites is sufficient, enabling operators to reduce deployment complexity while maintaining robust inspection capability.

\subsection{Modeling Assumptions \& Limitations}

This study assumes an offline planning problem with a uniform anomaly distribution along each pipeline. These assumptions are justified by the infrequency of failure events and the long planning horizons associated with subsea infrastructure deployment. Dynamic re-planning or multi-event scenarios are therefore not considered. The model can accommodate non-uniform failure distributions if such data are available by assigning higher weights to locations with higher failure likelihood, without changing the underlying optimization structure.

The results are naturally dependent on the geometry and density of the subsea pipeline system. We expect that configurations with elongated structures or geographically isolated pipeline systems will benefit more from the deployment of additional SDPs, whereas denser and more centrally connected systems will have a more rapid rate of saturation in terms of the maximum response time improvement. The configuration that was used in the study is a combination of a central hub and a peripheral pipeline system, similar to the Johan Sverdrup field. It is therefore representative of a larger, more mature subsea pipeline system. While the actual response time will vary from field to field, the trends observed in terms of diminishing returns for increasing fleet size and the importance of spatial placement for the worst-case scenario are likely to be applicable to similar subsea pipeline systems.

\subsection{Computational Scalability}
Table~\ref{tab:scalability} and Fig.~\ref{fig:scaling_curve} show that, even with 32‐core parallelism, the two‐phase MILP becomes impractical for $Q > 17$ and $P > 39$. At this point, average solve times exceed 10 minutes per $(R,S)$ scenario, with some low-resource cases (e.g., $R=1$, $S=2$) requiring over 80 minutes to solve individually. The number of binary variables grows superlinearly with problem size, reaching over 22 000 at $(Q=22, P=46)$. This results in significant computational load even for moderately sized oil and gas fields.

The bottleneck arises not only from the growing number of binary decision variables but also from the increased constraint coupling in configurations with small values of $R$, which are harder for the solver to resolve. These instances account for a disproportionate share of total runtime. While exhaustive MILP evaluation remains feasible for compact maps, scaling to larger environments or enabling interactive exploration requires algorithmic acceleration.

To restore tractability, several strategies can be explored. Warm starting from known previously solved $(R,S)$ configurations can reduce solver overhead. Heuristics can be used to eliminate dominated or clearly suboptimal configurations before invoking the full MILP. Alternatively, decomposition methods, surrogate models, or robust LP relaxations could be integrated while reducing solve time. These directions form the basis for future work toward scalable RAUV deployment planning.

\subsection{Future Work}
Building on the placement and assignment foundation established here, two primary directions are envisioned for future exploration.

Future work will examine scalability to larger subsea fields, denser anomaly discretizations, and multi-field infrastructure networks. Although the proposed MILP provides optimal solutions for the discretized instances considered here, larger cases may require decomposition, warm-starting, or related computational strategies to reduce solution time.

In practice, anomaly localization may be aided by the fusion of various pieces of information from different sources, which could include pressure monitoring systems, fiber optic sensing techniques like distributed acoustic sensing, acoustic monitoring techniques like hydrophone arrays, inspection histories, and operator reports. Such fusion may significantly reduce the search region for mission planning. This change may impact the problem description (Section~\ref{problem_statement}), as the current approach to the solution would require that spatially informed anomaly priors be incorporated. However, it is an open question how significantly this information would impact the optimal quantities and deployment locations of SDPs and RAUVs.

\section{CONCLUSION}
\label{conclusion}

We propose a two-phase MILP that jointly selects SDPs and RAUV home bases to minimize both maximum and average response times under anomaly location uncertainty. On a realistic, Johan Sverdrup-inspired network (and an extension of it with up to 17 candidate SDPs, 41 pipelines), the computational study uncovered three key insights: moderate redundancy (around 2 RAUVs and 4 SDPs) captured most of the gains in both metrics, a small set of central stations dominates optimal placement, and beyond the cost-time Pareto frontier further investment yields rapidly diminishing returns. 

Overall, the proposed framework offers offshore operators a rigorously benchmarked tool for balancing response time, asset investment in complex subsea environments.

\section{ACKNOWLEDGMENTS}

Funding for this work was provided by the NTNU VISTA Centre for Autonomous Robotic Operations Subsea (CAROS), a collaboration between Equinor, the Norwegian Academy of Science and Letters (DNVA), and the Norwegian University of Science and Technology (NTNU).

During the preparation of this work, the authors used Chat-GPT 4 and 5, and DeepL Write in order to improve language and readability. After using these tools/services, the authors reviewed and edited the content as needed and take full responsibility for the content of the publication.

\bibliographystyle{IEEEtran}
\bibliography{references}

\end{document}